%% file: Main.tex
\documentclass[12pt]{amsart}

\usepackage{amsmath}
\usepackage{amssymb}
\usepackage{amsfonts}
\usepackage{bbm}
\usepackage{mathrsfs}

\usepackage[usenames,dvipsnames]{color}
\usepackage{hyperref}
\usepackage{cite}

\newtheorem{definition}{Definition}

\newtheorem{theorem}{Theorem}

\theoremstyle{definition}
\newtheorem{remark}{Remark}

\newcommand{\eps}{\varepsilon}

\renewcommand{\theta}{\vartheta}
\newcommand{\Laplace}{\triangle}
\begin{document}
\title[]{Hyperuniform point sets on projective spaces}

\author{Johann S. Brauchart}
\address{Institute of Analysis and Number Theory,
  Graz University of Technology,
  Kopernikusgasse 24.
8010 Graz,
Austria}
\email{j.brauchart@tugraz.at}

\author{Peter J. Grabner\textsuperscript{\ddag}}
\thanks{\textsuperscript{\ddag} This author is supported by the Austrian
  Science Fund FWF project I~6750}
\email{peter.grabner@tugraz.at}
\maketitle
\begin{abstract}
  We extend the notion of hyperuniformity to the projective spaces
  $\mathbb{RP}^{d-1}$, $\mathbb{CP}^{d-1}$, $\mathbb{HP}^{d-1}$, and
  $\mathbb{OP}^2$. We show that hyperuniformity implies uniform distribution
  and present examples of deterministic point sets as well as point processes
  which exhibit hyperuniform behaviour.
\end{abstract}
\input{Introduction}
\input{Compact}
\input{Examples}

\bibliographystyle{amsplain}
\bibliography{refs}
\end{document}

%% file: Introduction.tex
\section{Introduction}\label{sec:introduction}

The quantitative characterisation of density fluctuation in many-particle systems led to the introduction of the concept of ``hyperuniformity'' (S. Torquato; cf. \cite{Torquato_Stillinger2003:local_density_fluctuations}) also called ``superhomogeneity'' (J.~Lebowitz, cf. \cite{Ghosh_Lebowitz2018:generalized_stealthy_hyperuniform_processes}).

% non-compact setting

The main feature of hyperuniformity of an infinite discrete point configuration $X$ in $\mathbb{R}^d$ is the fact that local density fluctuations are of smaller order than for a random (“Poissonian”) point configuration. In terms of the structure factor
\begin{equation*}
S(\mathbf{k})=\lim_{B\to\mathbb{R}^d}\frac1{\#(B\cap X)} \sum_{\mathbf{x},\mathbf{y}\in B\cap X}e^{i\langle\mathbf{k},\mathbf{x}-\mathbf{y}\rangle}\quad \text{(thermodynamic limit)}
\end{equation*}
this means that $\lim_{\mathbf{k}\to\mathbf{0}}S(\mathbf{k})=0$; i.e., normalised density fluctuations are completely suppressed at very large length-scales (cf.~\cite{Torquato_Stillinger2003:local_density_fluctuations}). This
\emph{thermodynamic limit} is understood in the sense that the volume $B$ (for instance a ball of radius $R$) tends to the whole space $\mathbb{R}^d$ while $\lim_{B\to\mathbb{R}^d}\frac{\#(B\cap  X)}{\mathrm{vol}(B)}=\rho$, where $\rho$ is the density. 
Equivalently, the number variance $\mathbb{V}[N_R]$ (understood in the sense of
thermodynamic limit) of particles within a spherical observation window of
radius $R$ can be used. When scaled by the window volume, this ratio tends to
zero as $R\to \infty$; i.e., $\mathbb{V}[N_R]$ grows more slowly than the
window volume of order $R^d$ in the large-$R$ limit. When
\begin{equation*}
S(\mathbf{k}) \sim \left\| \mathbf{k} \right\|^\alpha \qquad \text{as $\left\| \mathbf{k} \right\| \to \infty$,}
\end{equation*}
where $\alpha > 0$, the number variance has the following large-$R$ scaling: 
\begin{equation*}
\mathbb{V}[N_R] 
\sim
\begin{cases}
R^{d-1} & \alpha > 1, \\
R^{d-1} \log R & \alpha = 1, \\
R^{d-\alpha} & 0 < \alpha < 1. 
\end{cases}
\end{equation*}
This asymptotic behaviour provides a way to distinguish between different classes of hyperuniformity: Class I (strongest), II, and III (weakest). 
For a thorough discussion and survey regarding, in particular, hyperuniform states of matter, we refer to \cite{Torquato2018:hyperuniform_states_matter}. 
More recent contributions consider hyperuniformity in biology and geophysics \cite{Ge2023:hidden_order_turing_patterns,hu_Liu_Wa+2023:disordered,zhu2023dispersed}, cryptography \cite{gilpin2018:cryptographic}, and mathematics \cite{bjorklund_hartnick2023:hyperuniformity} to name a few.
% Coste, S.: Order, fluctuations, rigidities. https://scoste.fr/assets/survey_hyperuniformity.pdf (incomplete survey)
%
%
A nice layman's introduction to hyperuniformity can be found in~\cite{Wolchover2016:birds_eye_view_natures_hidden_order}.
%

% compact setting
In an answer to a question of S.~Torquato, we together with W.~Kusner~\cite{Brauchart_Grabner_Kusner2019:hyperuniform_deterministic} were able to extend the concept of hyperuniformity to the compact unit sphere $\mathbb{S}^d$ in the Euclidean space $\mathbb{R}^{d+1}$ by considering a fixed sequence of $N$-point distributions and analysing the asymptotic behaviour of the number variance with respect to certain spherical caps as $N \to \infty$. It turns out that there are three regimes of growth of the sperical caps as test sets that are of interest and only the one dubbed ``hyperuniform with respect to threshold order'' can be seen as a direct analogue to the classical notion of hyperuniformity; see  Definition~\ref{def:hyperuniformity.projective.space} below. 
In \cite{Brauchart_Grabner_Kusner2019:hyperuniform_deterministic} it is shown that hyperuniformity in each of the three regimes implies uniform distribution. Furthermore, it is proven that so-called Quasi-Monte Carlo design sequences (cf.~\cite{Brauchart_Saff_Sloan+2014:qmc-designs}) with strength at least $\frac{d+1}{2}$, certain energy minimising point set sequences (see, e.g., \cite{Brauchart_Grabner2015:distributing_sphere,Borodachov_Hardin_Saff2019:discreta_energy_on_rectifiable_sets}), and especially sequences of spherical designs of optimal growth order have hyperuniform behaviour. The paper~\cite{Brauchart_Grabner_Kusner+2020:hyperuniform_point_sets} studies hyperuniformity on the sphere for samples of point processes on the sphere. The spherical ensemble (cf.~\cite{Hough_Krishnapur_Peres+2009:zeros_gaussian}), the harmonic ensemble introduced in~\cite{Beltran_Marzo_Ortega-Cerda2016:determinantal}, and the jittered sampling process (shown to be a determinantal) exhibit hyperuniform behaviour. 

An essential requirement for our study of the hyperuniformity phenomenon in the compact setting is that the space is homogeneous. Having considered the sphere, it seems natural to turn to doubly homogeneous spaces next, which is the aim of this paper. We point out that another compact setting,  flat tori, is studied in~\cite{Stepanyuk2020:hyperuniform_point_sets}.

The paper is organised as follows. Section~\ref{sec:hyper-comp-doubly} introduces compact doubly homogeneous spaces, briefly discusses their harmonic analysis, and defines hyperuniformity on these spaces. Like in the spherical case \cite{Brauchart_Grabner_Kusner+2020:hyperuniform_point_sets, Brauchart_Grabner_Kusner2019:hyperuniform_deterministic}, there are three regimes. We arrive at the result that hyperuniformity, indeed, implies uniform distribution in each case. 
Section~\ref{sec:examples} discusses examples: maximisers of the sum of distances, jittered sampling, and the harmonic ensemble.

%%% Local Variables:
%%% mode: latex
%%% TeX-master: "Main"
%%% End:

%% file: Compact.tex
\section{Hyperuniformity on compact doubly homogeneous
  spaces}\label{sec:hyper-comp-doubly}
The main purpose of this paper is to extend the notion of hyperuniformity to
doubly homogeneous manifolds. Thus we first give a concise description of these
spaces and collect the relevant facts that will be used later.
\subsection{Harmonic analysis on compact doubly homogeneous
  spaces}
\label{sec:comp-doubly-homog}
Let us recall the definition of doubly homogeneous spaces. A metric space
$(X,\theta)$ is called doubly homogeneous, if there exists a group $G$ acting
isometrically on $X$ so that
\begin{equation*}
  \forall x_1,x_2,y_1,y_2\in X: \theta(x_1,x_2)=\theta(y_1,y_2)\Rightarrow
  \exists g\in G: y_1=gx_1\wedge y_2=gx_2.
\end{equation*}
The compact doubly homogeneous Riemannian manifolds have been characterised in
\cite{Wang1952:Two_Point_Compact} (see also
\cite{Helgason2000:Groups_geometric_analysis}):
\begin{align*}
  \mathbb{S}^d&\cong \mathrm{SO}(d+1)/\mathrm{SO}(d)\\
  \mathbb{RP}^{d-1}&\cong\mathrm{O}(d)/(\mathrm{O}(d-1)\times \mathrm{O}(1))\\
  \mathbb{CP}^{d-1}&\cong\mathrm{U}(d)/(\mathrm{U}(d-1)\times \mathrm{U}(1))\\
  \mathbb{HP}^{d-1}&\cong\mathrm{Sp}(d)/(\mathrm{Sp}(d-1)\times
  \mathrm{Sp}(1))\\
  \mathbb{OP}^2&\cong\mathrm{F}_4/\mathrm{Spin}(9).
\end{align*}
We also refer to \cite{Anderson_Dostert_Grabner+2023:riesz_green_energy} for
a concise introduction. The case of the sphere has been studied in
\cite{Brauchart_Grabner_Kusner+2020:hyperuniform_point_sets,
  Brauchart_Grabner_Kusner2019:hyperuniform_deterministic}; thus we mainly
focus on the projective spaces $\mathbb{FP}^{d-1}$ in this paper, where we use
$\mathbb{F}$ to represent the underlying scalar ``field'' ($\mathbb{R}$,
$\mathbb{C}$, $\mathbb{H}$, $\mathbb{O}$). The dimension of the space as a real
manifold is then given by
\begin{equation*}
  D=(d-1)\dim_{\mathbb{R}}(\mathbb{F}),
\end{equation*}
where
\begin{equation*}
  \dim_{\mathbb{R}}(\mathbb{C})=2,\quad
  \dim_{\mathbb{R}}(\mathbb{H})=4,\quad \dim_{\mathbb{R}}(\mathbb{O})=8. 
\end{equation*}
Furthermore, we associate the parameters
\begin{equation*}
  \alpha=\frac{d-1}2\dim_{\mathbb{R}}(\mathbb{F})-1\quad\text{and }
  \beta=\frac12\dim_{\mathbb{R}}(\mathbb{F})-1
\end{equation*}
to the spaces $\mathbb{FP}^{d-1}$. We normalise the geodesic metric $\theta$ to
take values in $[0,\frac\pi2]$, and equip the space with the normalised surface
measure $\sigma$.

The measure of the geodesic ball
$B(x,r)=\{y\in\mathbb{FP}^{d-1}\mid\theta(x,y)<r\}$ is then given by
\begin{equation*}
  \sigma(B(x,r))=\int_0^r A(\theta)\,d\theta,
\end{equation*}
where
\begin{equation*}
  A(r)=\frac{2\Gamma(\alpha+\beta+2)}{\Gamma(\alpha+1)\Gamma(\beta+1)}
  \sin(r)^{2\alpha+1}\cos(r)^{2\beta+1}
\end{equation*}
denotes the surface area of the geodesic sphere
$\mathbb{S}(a,r)=\{x\in\mathbb{FP}^{d-1}\mid \theta(a,x)=r\}$. From now on, we
use the notation
\begin{equation*}
  C_{\alpha,\beta}=
  \frac{2\Gamma(\alpha+\beta+2)}{\Gamma(\alpha+1)\Gamma(\beta+1)}
\end{equation*}
for the normalising constant.

The Laplace operator associated to the underlying metric tensor is then given
by
\begin{equation*}
  \Laplace f=-\frac1{A(r)}\frac\partial{\partial r}
  \left(A(r)\frac{\partial f}{\partial r}\right)+\Laplace_{\mathbb{S}(a,r)}f,
\end{equation*}
where $\Laplace_{\mathbb{S}(a,r)}$ denotes the Laplace operator on the geodesic
sphere. In the following we will restrict our attention to zonal functions
depending only on $r$; for these functions $\Laplace_{\mathbb{S}(a,r)}f$
vanishes. Furthermore, we notice that we use the geometers convention that the
Laplacian has non-negative eigenvalues.

The zonal eigenfunctions of $\Laplace$ are then given by the functions
\begin{equation*}
  P_n^{(\alpha,\beta)}(\cos(2\theta(x,a))),
\end{equation*}
the Jacobi-polynomials with parameters $\alpha$ and $\beta$. The corresponding
eigenvalues are given by
\begin{equation*}
  \lambda_n=4n(n+\alpha+\beta+1).
\end{equation*}
The dimension of the corresponding space of eigenfunctions $V_n$ is then given
by
\begin{equation*}
  m_n=\frac{2n+\alpha+\beta+1}{\alpha+\beta+1}
  \frac{(\alpha+\beta+1)_n(\alpha+1)_n}{n!(\beta+1)_n},
\end{equation*}
where $(x)_n=x(x+1)\cdots(x+n-1)$ denotes the rising factorial (Pochhammer
symbol).

\begin{table}[h]
  \renewcommand{\arraystretch}{1.8}
  \begin{tabular}{l| l| l|l|l}
    $X$  & $\alpha$&$\beta$& $\lambda_k$&$m_k = \dim(V_k)$ \\[2mm]\hline
    $\mathbb{RP}^{d-1}$ & $\frac{d-3}2$&$-\frac12$&$2k (2k+d-2)$ &
    $ \frac{4k+d-2}{d-2}\binom{2k+d-3}{d-3}$\\[2mm]\hline
    $\mathbb{CP}^{d-1}$ & $d-2$&$0$&$4k (k+d-1)$ &
    $ \frac{2k+d-1}{d-1}\binom{d+k-2}{d-2}^2$\\[2mm]\hline
    $\mathbb{HP}^{d-1}$&$2d-3$&$1$ & $4k (k+2d-1)$ &
    $ \frac{2k+2d-1}{(2d-1)(2d-2)} \binom{k+2d-2}{2d-2} \binom{k+2d-3}{2d-3}$
    \\[2mm]\hline
    $\mathbb{OP}^{2}$ & $7$&$3$&$4k (k+11)$ &
    $ \frac{2k+11}{1320} \binom{k+10}{7} \binom{k+7}{7}$\\
  \end{tabular}
  \caption{\label{tab:eigen}
    The eigenvalues and dimensions of the eigenspaces of the
    Laplace operator on the projective spaces $\mathbb{FP}^{d-1}$}
\end{table}

The spaces of eigenfunctions $V_n$ are invariant under the group $G$ acting on
$\mathbb{FP}^{d-1}$. The corresponding representation turns out to be
irreducible. Let $Y_{n,k}(x)$ ($k=1,\ldots,m_n$) denote an orthonormal basis of
$V_n$. Then the addition theorem
\begin{equation}\label{eq:addition}
  \sum_{k=1}^{m_n}Y_{n,k}(x)Y_{n,k}(y)=\frac{m_n}{P_n^{(\alpha,\beta)}(1)}
  P_n^{(\alpha,\beta)}(\cos(2\theta(x,y)))
\end{equation}
is an immediate consequence of the irreducibility of $V_n$. Furthermore, we
recall the formula for integrals of zonal functions
\begin{equation}
  \label{eq:zonal}
  \begin{split}
    &\int_{\mathbb{FP}^{d-1}}f(\cos(2\theta(x,y)))\,d\sigma(x)\\
    =
    &C_{\alpha,\beta}\int_0^{\frac\pi2}f(\cos(2\theta))\sin(\theta)^{2\alpha+1}
    \cos(\theta)^{2\beta+1}\,d\theta,
  \end{split}
\end{equation}
see for instance \cite{Helgason1965:radon_transform}.

Every function 
$f(\cos(2t))\in L^2([0,\frac\pi2],\sin(t)^{2\alpha+1}\cos(t)^{2\beta+1}\,dt)$
can be expanded in terms of the orthogonal system
$P_n^{(\alpha,\beta)}(\cos(2t))$:
\begin{equation}\label{eq:fourier}
  f(\cos(2t))=\sum_{n=0}^\infty \widehat{f}(n)P_n^{(\alpha,\beta)}(\cos(2t)),
\end{equation}
where the coefficients are given by
\begin{equation*}
  \widehat{f}(n)=\frac{C_{\alpha,\beta}m_n}{\left(P_n^{(\alpha,\beta)}(1)\right)^2}
  \int\limits_0^{\frac\pi 2}f(\cos(2t))P_n^{(\alpha,\beta)}(\cos(2t))
  \sin(t)^{2\alpha+1}\cos(t)^{2\beta+1}\,dt.
\end{equation*}
The convergence in \eqref{eq:fourier} is \emph{a priori} only in the
$L^2$-sense; for continuous functions $f$ satisfying
\begin{equation}\label{eq:pos-def}
  \sum_{j,k=1}^Nc_j\overline{c_k}f(\cos(2\theta(x_j,x_k)))\geq0
\end{equation}
for all choices of $N$, $x_1,\ldots,x_N\in X$, and
$c_1,\ldots,c_N\in\mathbb{C}$, Mercer's theorem ensures absolute and uniform
convergence of the series \eqref{eq:fourier} (see for instance
\cite{Ferreira_Menegatto2009:eigenvalues_integral_operators}). Functions
satisfying \eqref{eq:pos-def} are called \emph{positive definite}. It follows
from \cite{Bochner1941:Positive_Definite} that positive definite functions are
characterised by the non-negativity of the coefficients $\widehat{f}(n)$.

\subsection{Hyperuniformity on $\mathbb{FP}^{d-1}$}
We use similar ideas as in
\cite{Brauchart_Grabner_Kusner+2020:hyperuniform_point_sets,
  Brauchart_Grabner_Kusner2019:hyperuniform_deterministic}
to define hyperuniform sequences of point sets in the setting of the spaces
$\mathbb{FP}^{d-1}$. As in the original euclidean setting the concept is based
on the \emph{number variance}.
\begin{definition}\label{def-hyper}
  Let $(X_N)_{N\in\mathbb{N}}$ a sequence of point sets in $\mathbb{FP}^{d-1}$
  with $\lim_{N\to\infty}\#X_N=\infty$. The number variance of this sequence is
  then given by
  \begin{equation}
    \label{eq:variance}
    \begin{split}
      &\mathbb{V}(X_N,r)=\mathbb{V}_x\#(X_N\cap B(x,r))\\
      =&
      \int_{\mathbb{FP}^{d-1}}\left(\sum_{y\in X_N}\mathbbm{1}_{B(x,r)}(y)-
        \#X_N\sigma(B(x,r))\right)^2\,d\sigma(x).
    \end{split}
  \end{equation}
\end{definition}
In a probabilistic setting the number variance is the variance of the number of
points in a geodesic ball of radius $r$. The interest in this quantity stems
from the fact that it can be used as a measure for the quality of a point
distribution in the deterministic as well as the probabilistic setting.

\begin{remark}
  The number variance is intimately related to the $L^2$-discrepancy
  \begin{equation*}
    \mathrm{D}_{L^2}(X_N)=\frac1{\#X_N}
    \left(\int_0^{\frac\pi2}\mathbb{V}(X_N,r)\,dr\right)^{\frac12}.
  \end{equation*}
  This notion is a classical measure for the deviation of the discrete
  distribution $\frac1{\#X_N}\sum_{x\in X_N}\delta_x$ from uniform
  distribution. $L^2$-discrepancy has been studied to find lower bounds for this
  deviation; the theory of \emph{irregularities of distribution} has developed
  techniques to derive lower bounds for various notions of discrepancy. For a
  comprehensive introduction and a very good overview over results in this
  context we refer to \cite{Beck_Chen2008:irregularities_distribution}.
\end{remark}
\begin{remark}\label{rem:invariance}
  Stolarsky's celebrated invariance principle
  \cite{Stolarsky1973:sums_distances_sphere} relates the $L^2$-discrepancy of
  a point set $X_N$ on the sphere $\mathbb{S}^d$ to the sum of mutual chordal
  distances
  \begin{equation}\label{eq:invariance}
    \mathrm{D}_{L^2}(X_N)^2+c_d\frac1{(\#X_N)^2}\sum_{x,y\in X_N}\|x-y\|=
    c_d\iint_{\mathbb{S}^d\times\mathbb{S}^d}\|x-y\|\,d\sigma(x)\,d\sigma(y)
  \end{equation}
  for a suitable constant $c_d$. This shows that the $L^2$-discrepancy is
  minimised for sets maximising the sum of mutual distances.

  This invariance principle was generalised to the projective spaces
  $\mathbb{FP}^{d-1}$ by Skriganov
  \cite{Skriganov2020:stolarskys_invariance_principle}. Furthermore, it is
  observed there that the invariance principle is a special case of a much more
  general relation.
\end{remark}

As in
the classical notion of hyperuniformity we will be interested in the dependence
of the number variance on the radius $r$. Three regimes turn out to be of
interest.
\begin{definition} \label{def:hyperuniformity.projective.space}
  A sequence of point sets $(X_N)_{N\in\mathbb{N}}$ is called
  \begin{enumerate}
  \item hyperuniform for large balls, if for every $r\in(0,\frac\pi2)$
    \begin{equation*}
      \mathbb{V}(X_N,r)=o\left(\#X_N\right),
    \end{equation*}
  \item hyperuniform for small balls, if for every sequence of
    $(r_N)_{N\in\mathbb{N}}$ of radii satisfying
    \begin{align*}
      &\lim_{N\to\infty}r_N=0\\
      &\lim_{N\to\infty}\sigma(B(\cdot,r_N))\#X_N=\infty
    \end{align*}
    the following holds
    \begin{equation*}
      \mathbb{V}(X_N,r)=o\left(\#X_N\sigma(B(\cdot,r_N))\right),
    \end{equation*}
  \item hyperuniform for balls at threshold order, if
    \begin{equation*}
      \limsup_{N\to\infty}\mathbb{V}(X_N,r(\#X_N)^{-\frac1D})=\mathcal{O}(r^{D-1})
    \end{equation*}
    for $r\to\infty$.
  \end{enumerate}
\end{definition}
We notice that hyperuniformity at threshold order is obtained by rescaling and
adapting the classical notion to the compact situation. The motivation for the
order $r^{D-1}$ on the right hand side comes from the order of the surface area
of the geodesic sphere; i.i.d. random points would achieve the order
$r^D$. Furthermore, it is clear that a smaller order of magnitude than
$(\#X_N)^{-\frac1D}$ for the radii would not be sensible, because then the
expected number of points in a ball would tend to $0$.  The other two notions
of hyperuniformity have no obvious counterpart in the classical setting. All
three notions have in common that the number variance is required to have
smaller order of magnitude than the variance of i.i.d. random point sets of the
same cardinality.
\subsection{Uniform distribution}\label{sec:uniform-distribution}
Uniform distribution is a notion that formalises the discretisation of
measures.
\begin{definition}
  A sequence of of point sets $(X_N)_{N\in\mathbb{N}}$ in $\mathbb{FP}^{d-1}$
  with $\lim_{N\to\infty}\#X_N=\infty$ is called uniformly distributed, if
  \begin{equation*}
    \lim_{N\to\infty}\frac1{\#X_N}\sum_{x\in
      X_N}\mathbbm{1}_{B(x,r)}=\sigma(B(x,r))
  \end{equation*}
  for all $x\in \mathbb{FP}^{d-1}$ and $r\in[0,\frac\pi2]$.
\end{definition}
For a comprehensive introduction to uniform distribution we refer to
\cite{Drmota_Tichy1997:sequences_discrepancies_applications,
  Kuipers_Niederreiter1974:uniform_distribution_sequences}.

From the general theory (see
\cite{Kuipers_Niederreiter1974:uniform_distribution_sequences}) it is known
that a sequence of point sets is uniformly distributed, if and only if
\begin{equation*}
  \lim_{N\to\infty}\frac1{\#X_N}\sum_{x\in X_N}f(x)=
  \int_{\mathbb{FP}^{d-1}}f(x)\,d\sigma(x)
\end{equation*}
for all continuous functions $f:\mathbb{FP}^{d-1}\to\mathbb{C}$. From this it
follows immediately (using the denseness of the Laplace eigenfunctions and the
addition theorem \eqref{eq:addition}) that $(X_N)_{N\in\mathbb{N}}$ is
uniformly distributed, if and only if
\begin{equation}\label{eq:Weyl}
  \lim_{N\to\infty}\frac1{(\#X_N)^2}\sum_{x,y\in X_N}P_n^{(\alpha,\beta)}
  (\cos(2\theta(x,y)))=0
\end{equation}
for all $n\geq1$.

We now want to relate hyperuniformity to uniform distribution. Since uniform
distribution is a property modelled after the law of large numbers, which holds
especially for i.i.d. random points, and hyperuniformity is modelled as a
property ``better than i.i.d. random points'', we can only expect that
hyperuniformity implies uniform distribution.

In order to show that in all three regimes, we first derive an expression for
the number variance of $X_N$ in terms of Jacobi polynomials. We observe that
$\mathbbm{1}_{B(x,r)}(y)=\mathbbm{1}_{[0,r)}(\theta(x,y))$ is a zonal
function. Thus we can expand it as a series in terms of Jacobi polynomials in
$\cos(2\theta(x,y))$. Using the formula
\begin{multline*}
  \int_0^rP_n^{(\alpha,\beta)}(\cos(2t))\sin(t)^{2\alpha+1}\cos(t)^{2\beta+1}\,dt\\
  =  \frac1{2n}\sin(r)^{2\alpha+2}\cos(r)^{2\beta+2}
  P_{n-1}^{(\alpha+1,\beta+1)}(\cos(2r))
\end{multline*}
for $n\geq1$ (see
\cite{Magnus_Oberhettinger_Soni1966:formulas_theorems_special}), we obtain
\begin{equation*}
  \mathbbm{1}_{B(x,r)}(y)=\sigma(B(x,r))+
  \sum_{n=1}^\infty\frac{C_{\alpha,\beta}m_n}{P_n^{\alpha,\beta}(1)}a_n(r)
  P_n^{(\alpha,\beta)}(\cos(2\theta(x,y))).
\end{equation*}
with
\begin{equation*}
  a_n(r)=\frac1{2nP_n^{\alpha,\beta}(1)}\sin(r)^{2\alpha+2}\cos(r)^{2\beta+2}
  P_{n-1}^{(\alpha+1,\beta+1)}(\cos(2r)).
\end{equation*}

This equation is valid in the $L^2$-sense. Applying Parseval's identity we
obtain
\begin{equation}
  \label{eq:var-Jacobi}
\mathbb{V}(X_n,r)=
\sum_{n=1}^\infty
    \frac{C_{\alpha,\beta}^2m_n}{P_n^{\alpha,\beta}(1)}a_n(r)^2
    \sum_{x,y\in X_N}P_n^{(\alpha,\beta)}(\cos(2\theta(x,y))).
\end{equation}
This series is absolutely and uniformly (as a function of $x$ and $y$)
convergent by Mercer's theorem. Notice that the function
\begin{equation}\label{eq:g_r-formula}
  \begin{split}
    g_r(x,y)&=\!\!\!\!\int\limits_{\mathbb{FP}^{d-1}}\!\!\!
    \left(\mathbbm{1}_{B(x,r)}(z)-\sigma(B(x,r))\right)
    \left(\mathbbm{1}_{B(y,r)}(z)-\sigma(B(y,r))\right)\,d\sigma(z)\\
    &=\sigma(B(x,r)\cap B(y,r))-\sigma(B(\cdot,r))^2
  \end{split}
\end{equation}
is positive definite and that
\begin{equation}\label{eq:v=gr}
  \mathbb{V}(X_N,r)=\sum_{x,y\in X_N}g_r(x,y).
\end{equation}

\subsubsection{Hyperuniformity for large balls implies uniform
  distribution}\label{sec:hyper-large-balls}
Assume now that $(X_N)_{N\in\mathbb{N}}$ is hyperuniform for large balls. For
fixed $n$ choose a distance $0<r<\frac\pi2$ such that
$P_n^{(\alpha,\beta)}(\cos(2r))\neq0$. Then equation \eqref{eq:var-Jacobi}
yields
\begin{multline*}
  0=\lim_{N\to\infty}\frac{\mathbb{V}(X_N,r)}{\#X_N}\geq\\
  \frac{C_{\alpha,\beta}^2m_n}{P_n^{\alpha,\beta}(1)}a_n(r)^2
  \lim_{N\to\infty}\frac1{\#X_N}
  \sum_{x,y\in X_N}P_n^{(\alpha,\beta)}(\cos(2\theta(x,y))).
\end{multline*}
This proves
\begin{theorem}\label{thm:large}
  Let $(X_N)_{N\in\mathbb{N}}$ be hyperuniform for large balls. Then
  $(X_N)_{N\in\mathbb{N}}$ is uniformly distributed. More precisely,
  \begin{equation}\label{eq:Weyl-better}
    \lim_{N\to\infty}\frac1{\#X_N}
  \sum_{x,y\in X_N}P_n^{(\alpha,\beta)}(\cos(2\theta(x,y)))=0
\end{equation}
holds for all $n\geq1$.
\end{theorem}
\begin{remark}
  Notice that \eqref{eq:Weyl-better} gives a better rate of convergence as
  compared to the criterion \eqref{eq:Weyl} for uniform distribution. This
  reflects the heuristic expectation that hyperuniformity should give an
  improved quality of uniform distribution.
\end{remark}

\subsubsection{Hyperuniformity for small balls implies uniform
  distribution}\label{sec:hyper-small-balls}
Assume that $(X_N)_{N\in\mathbb{N}}$ is hyperuniform for small
balls. Furthermore, we observe that for fixed $n\geq1$ the coeffient of
$P_n^{(\alpha,\beta)}$ in \eqref{eq:var-Jacobi} is of order $r_N^{2D}$ for
$N\to\infty$ (notice that $4\alpha+4=2D$). Arguing as above we obtain that
\begin{equation*}
  \lim_{N\to\infty}\frac{\mathbb{V}(X_N,r_N)}{\#X_N\sigma(B(\cdot,r_N))}=0
\end{equation*}
implies
\begin{equation*}
  \lim_{N\to\infty}\frac{r_N^D}{\#X_N}
  \sum_{x,y\in X_N}P_n^{(\alpha,\beta)}(\cos(2\theta(x,y)))=0.
\end{equation*}
Since $r_N$ can be chosen to tend to $0$ arbitrarily slowly, this implies
\begin{equation*}
  \limsup_{N\to\infty}\frac{1}{\#X_N}
  \sum_{x,y\in X_N}P_n^{(\alpha,\beta)}(\cos(2\theta(x,y)))<\infty.
\end{equation*}
Thus we have shown
\begin{theorem}\label{thm:small}
  Let $(X_N)_{N\in\mathbb{N}}$ be hyperuniform for small balls. Then
  $(X_N)_{N\in\mathbb{N}}$ is uniformly distributed. More precisely,
  \begin{equation}\label{eq:Weyl-better2}
    \limsup_{N\to\infty}\frac1{\#X_N}
  \sum_{x,y\in X_N}P_n^{(\alpha,\beta)}(\cos(2\theta(x,y)))<\infty
\end{equation}
holds for all $n\geq1$.
\end{theorem}
\begin{remark}
  Notice that we have an improvement in the quality of uniform distribution in
  comparison with \eqref{eq:Weyl} by a power $(\#X_N)^{1-\eps}$ for $\eps>0$.
\end{remark}
\subsubsection{Hyperuniformity for balls at threshold order implies
  uniform distribution}\label{sec:hyper-balls-threshold}
Assume that $(X_N)_{N\in\mathbb{N}}$ is hyperuniform for balls at threshold
order. For fixed $n\geq 1$ we notice that
\begin{equation*}
  a_n(r(\#X_N)^{-\frac1D})^2
  \sim
  \frac{P_{n-1}^{(\alpha+1,\beta+1)}(1)^2}{(2nP_n^{(\alpha,\beta)}(1))^2}
  \frac{r^{2D}}{(\#X_N)^2}.
\end{equation*}
Combining this with \eqref{eq:var-Jacobi} gives
\begin{multline*}
  r^{2D}\frac{C_{\alpha,\beta}^2m_nP_{n-1}^{(\alpha+1,\beta+1)}(1)^2}
  {(2n)^2(P_n^{\alpha,\beta}(1))^3}
  \limsup_{N\to\infty}\frac1{(\#X_N)^2}
  \sum_{x,y\in X_N}P_n^{(\alpha,\beta)}(\cos(2\theta(x,y)))\\
  \leq
  \limsup_{N\to\infty}\mathbb{V}(X_N,r(\#X_N)^{-\frac1D})=\mathcal{O}(r^{D-1}).
\end{multline*}
This can only hold, if
\begin{equation*}
  \limsup_{N\to\infty}\frac1{(\#X_N)^2}
  \sum_{x,y\in X_N}P_n^{(\alpha,\beta)}(\cos(2\theta(x,y)))=0.
\end{equation*}
Summing up, we have proved
\begin{theorem}
  Let $(X_N)_{N\in\mathbb{N}}$ be hyperuniform for balls at threshold
  order. Then $(X_N)_{N\in\mathbb{N}}$ is uniformly distributed.
\end{theorem}
\begin{remark}
  In this case there seems to be no improvement in the order of convergence of
  the sum \eqref{eq:Weyl}.
\end{remark}
%%% Local Variables:
%%% mode: latex
%%% TeX-master: "Main"
%%% End:

%% file: Examples.tex
\section{Examples}\label{sec:examples}

\subsection{Maximisers of the sum of distances}
\label{sec:maxim-sum-dist}
As pointed out in Remark~\ref{rem:invariance} Skriganov
\cite{Skriganov2020:stolarskys_invariance_principle} generalised Stolarsky's
invariance principle \cite{Stolarsky1973:sums_distances_sphere} to the
projective spaces $\mathbb{FP}^{d-1}$. As a consequence the $L^2$-discrepancy
is minimised amongst all sets of the same cardinality by configurations
maximising the sum of mutual ``chordal'' distances
\begin{equation}\label{eq:sum-dist}
  \sum_{x,y\in X_N}\sin(\theta(x,y)).
\end{equation}
It was shown in \cite{Skriganov2019:point_distributions_homogeneous} that there
exist positive constants $K_d$ and $k_d$ such that
\begin{equation}
  \label{eq:skriganov}
  \begin{split}
  -&K_d(\#X_N)^{1-\frac1D}\leq\\
  &\sum_{x,y\in X_N}\sin(\theta(x,y))-(\#X_N)^2\iint\limits_{(\mathbb{FP}^{d-1})^2}
  \sin(\theta(x,y))\,d\sigma(x)\,d\sigma(y)\\
  &\leq-k_d(\#X_N)^{1-\frac1D}
  \end{split}
\end{equation}
for a point set $X_N$ maximising \eqref{eq:sum-dist}.

Integrating \eqref{eq:var-Jacobi} against $\sin(2r)\,dr$ and using
\begin{align*}
  &\int_0^{\frac\pi2}\sin(r)^{4\alpha+4}\cos(r)^{4\beta+4}
  P_{n-1}^{(\alpha+1,\beta+1)}(\cos(2r))^2\sin(2r)\,dr=\\
  &\frac{\Gamma(2\alpha+3)\Gamma(2\beta+3)}{\alpha+\beta+2}
  2^{2n-2}(n+\alpha+\beta+1)\times\\
  &\frac{(\frac12)_{n-1}(\alpha+2)_{n-1}(\beta+2)_{n-1}(\alpha+\beta+2)_{n-1}}
  {((n-1)!)^2\Gamma(2(n+\alpha+\beta+3))}
\end{align*}
(see \cite{Skriganov2020:stolarskys_invariance_principle})
we obtain
\begin{equation}
  \label{eq:l2-discr}
  \begin{split}
    &\mathrm{D}_{L^2}(X_N)^2=C_{\alpha,\beta}^2\Gamma(2\alpha+3)\Gamma(2\beta+3)
    \times\\
    &\sum_{n=1}^\infty 2^{2n-2}(2n+\alpha+\beta+1)
    \frac{(\frac12)_{n-1}(\alpha+\beta+1)_n(\alpha+\beta+2)_n}
    {(\alpha+1)_n\Gamma(2(n+\alpha+\beta+3))}\times\\
    &\sum_{x,y\in X_N}
    P_n^{(\alpha,\beta)}(\cos(2\theta(x,y)))
  \end{split}
\end{equation}
(see also \cite{Skriganov2020:stolarskys_invariance_principle,
  Skriganov2019:stolarskys_invariance_principle_ii}).
We notice that the coefficients in the sum are of the order $n^{-\alpha-2}$.

We will use the inequality
\begin{equation}\label{eq:jacobi-bound}
  \begin{split}
    &\sin(\theta)^{a+\frac12}\cos(\theta)^{b+\frac12}
    |P_n^{(a,b)}(\cos(2\theta))|\\
    &\leq\frac{\Gamma(a+1)}{\sqrt\pi}
    \binom{n+a}n
    \left(n+\frac12(a+b+1)\right)^{-a-\frac12}
    =\mathcal{O}\left( n^{-\frac12}\right)
  \end{split}
\end{equation}
valid for $a\geq b>-1$ (see
\cite{Chow_Gatteschi_Wong1994:bernstein_type_inequality}). Inserting
$a=\alpha+1$ and $b=\beta+1$ in \eqref{eq:jacobi-bound} we estimate the number
variance by
\begin{equation}\label{eq:variance-upper}
  \mathbb{V}(X_N,r)\leq C\sin(r)^{2\alpha+1}\cos(r)^{2\beta+1}
  \sum_{n=1}^\infty n^{-\alpha-2}
  \sum_{x,y\in X_N}P_n^{(\alpha,\beta)}(\cos(2\theta(x,y))).
\end{equation}
Putting \eqref{eq:l2-discr} and~\eqref{eq:variance-upper} together and
observing that $2\alpha+1=D-1$ we obtain
\begin{equation}
  \label{eq:hyper}
  \mathbb{V}(X_N,r)=\mathcal{O}\left( \sin(r)^{D-1}\mathrm{D}_{L^2}(X_N)\right).
\end{equation}
For point sets with optimal $L^2$-discrepancy (or equivalently, maximal sum of
distances) we obtain
\begin{equation*}
  \mathbb{V}(X_N,r)=\mathcal{O}\left( \sin(r)^{D-1}(\#X_N)^{1-\frac1D}\right).
\end{equation*}
This inequality immediately implies the theorem.
\begin{theorem}
  Let $(X_N)_{N\in\mathbb{N}}$ be a sequence of maximisers of the sum of mutual
  distances \eqref{eq:sum-dist}. Then this sequence is hyperuniform in all all
  three regimes.
\end{theorem}
\subsection{Jittered sampling}\label{sec:jittered-sampling}
A very simple and obvious probabilistic model for a point set $X_N$ is given by
jittered sampling. We start with a partition of the space $\mathbb{FP}^{d-1}$
into $N$ sets $A_j$ ($j=1,\ldots,N$) of equal area and diameter
$\leq CN^{-\frac1D}$, where the constant only depends on the space. It is known
from \cite{Gigante_Leopardi2017:diameter_ahlfors} that such a partition exists.

We define a point process $\mathscr{X}_N^J$, called the \emph{jittered sampling
  process}, by choosing $N$ points $x_j\in A_j$ ($j=1,\ldots,N$), where each
point is chosen with respect to the surface measure restricted to $A_j$. The
number variance of this process is then given by
\begin{equation}\label{eq:variance-jittered}
  \begin{split}
    &V_N=\mathbb{V}(\mathscr{X}_N^J,r)\\
    &=\!\!\! \int\limits_{\mathbb{FP}^{d-1}}\int\limits_{A_1}\!\!\cdots\!\!\int\limits_{A_N}
    \!\!\left(\sum_{j=1}^N\mathbbm{1}_{B(x,r)}(x_j)-N\sigma(B(x,r))\right)^2\!\!
    d\sigma_1(x_1)\cdots d\sigma_N(x_N)\,d\sigma(x),
  \end{split}
\end{equation}
where $\sigma_j(A)=N\sigma(A_j\cap A)$ denotes the surface measure restricted
to $A_j$.

Expanding the square in~\eqref{eq:variance-jittered} and using the fact that
$\mathbbm{1}_{B(y,r)}(x)=\mathbbm{1}_{B(x,r)}(y)$ we obtain
\begin{align*}
  V_N&=\sum_{i\neq j}\int_{A_i}\int_{A_j}
  \sigma(B(x_i,r)\cap B(x_j,r))\,d\sigma_i(x_i)\,
  d\sigma_j(x_j)\\
  &+N\sigma(B(\cdot,r))-N^2\sigma(B(\cdot,r))^2 \\
  &=N^2\int_{\mathbb{FP}^{d-1}}\int_{\mathbb{FP}^{d-1}}
  \sigma(B(x,r)\cap B(y,r))\,d\sigma(x)\,
  d\sigma(y)\\
  &-
  \sum_{i=1}^N\int_{A_i}\int_{A_i}\sigma(B(x_i,r)\cap B(x_j,r))\,d\sigma_i(x_i)\,
  d\sigma_j(x_j)\\
  &+N\sigma(B(\cdot,r))-N^2\sigma(B(\cdot,r))^2\\
  &=\frac12\sum_{i=1}^N\int_{A_i}\int_{A_i}\sigma(B(x_i,r)\triangle B(y_i,r))
  \,d\sigma_i(x_i)\,d\sigma_i(y_i),
\end{align*}
where we have used
\begin{equation*}
  \int_{\mathbb{FP}^{d-1}}\int_{\mathbb{FP}^{d-1}}
  \sigma(B(x,r)\cap B(y,r))\,d\sigma(x)\,
  d\sigma(y)=\sigma(B(\cdot,r))^2
\end{equation*}
and used the notation $A\triangle B=A\cup B\setminus(A\cap B)$ for the
symmetric difference. A similar computation has been used in
\cite{Brauchart_Grabner_Kusner+2020:hyperuniform_point_sets}.

The measure of the symmetric difference of two balls can be bounded
\begin{equation}\label{eq:symmetric}
  \sigma(B(x,r)\triangle B(y,r))\leq
  C'\theta(x,y)\mathrm{surface}(\partial B(\cdot,r))
\end{equation}
for some constant $C'>0$. Since the diameter of each part $A_j$ of the partition
is bounded by $CN^{-\frac1D}$ we obtain the bound
\begin{equation*}
  V_N\leq CC'N^{1-\frac1D}\mathrm{surface}(\partial B(\cdot,r))=
  \mathcal{O}(r^{D-1}N^{1-\frac1D}),
\end{equation*}
from which we derive the theorem.
\begin{theorem}
  The jittered sampling process $\mathscr{X}_N^J$ is hyperuniform in all
three regimes.
\end{theorem}
\subsection{The harmonic ensemble}\label{sec:harmonic-ensemble}
Determinantal point processes (see
\cite{Hough_Krishnapur_Peres+2009:zeros_gaussian}) have become a very useful
tool to provide probabilistic models, which behave better than i.i.d. models
due to built mutual repulsion of the points. Especially, in the context of
finding well distributed point sets such processes have proven to be a very
useful tool to establish good bounds for various quality measures of point
distributions (see, for instance
\cite{Anderson_Dostert_Grabner+2023:riesz_green_energy,
  Beltran_Marzo_Ortega-Cerda2016:determinantal,
  Beltran-Etayo2018:projective_ensemble, Hirao2021:finite_frames_determinantal,
  Beltran_Ferizovic2020:approximation_uniform_SO3,
  Marzo_Ortega-Cerda2018:expected_riesz_energy,
  Hirao2018:qmc_designs_determinantal}). Since the notion of hyperuniformity
has been motivated as a property of point sets, which behave ``better than
i.i.d. random'', it is natural to test samples of determinantal processes for
this property.

The rotation invariant processes on the spaces $\mathbb{FP}^{d-1}$ have been
characterised in \cite{Anderson_Dostert_Grabner+2023:riesz_green_energy}. The
most natural amongst them is the process induced by the projection kernel on
the spaces of eigenfunctions for the $N+1$ smallest eigenvalues of the Laplace
operator. This kernel is given by
\begin{equation}
  \label{eq:harmonic-kernel}
  \begin{split}
    K_N(x,y)&=\sum_{n=0}^N\sum_{k=1}^{m_n}Y_{n,k}(x)Y_{n,k}(y)=
    \sum_{n=0}^N\frac{m_n}{P_n^{(\alpha,\beta)}(1)}
    P_n^{(\alpha,\beta)}(\cos(2\theta(x,y)))\\
    &=
    \frac{(\alpha+\beta+2)_N}{(\beta+1)_N}P_N^{(\alpha+1,\beta)}(\cos(2\theta(x,y)))
  \end{split}
\end{equation}
(see, for instance \cite{Andrews_Askey_Roy1999:special_functions}). This
process produces
\begin{equation*}
  K_N(x,x)=\frac{(\alpha+\beta+2)_N(\alpha+2)_N}{N!(\beta+1)_N}\sim
  \frac{\Gamma(\beta+1)}{\Gamma(\alpha+\beta+2)\Gamma(\alpha+2)}N^{2\alpha+2}
\end{equation*}
points; we use the notation $\mathscr{X}_N^H$ for this process and call it the
\emph{harmonic process} like the corresponding process on the sphere introduced
in \cite{Beltran_Marzo_Ortega-Cerda2016:determinantal}. For more details
on determinantal point processes we refer to
\cite{Hough_Krishnapur_Peres+2009:zeros_gaussian} again.

We notice that the number variance of $\mathscr{X}_N^H$ equals the expectation
of the sum \eqref{eq:v=gr} under law of the process $\mathscr{X}_N^H$,
\begin{equation}
  \label{eq:var-harmonic}
  \mathbb{V}(\mathscr{X}_N^H,r)=\mathbb{E}\sum_{x,y}g_r(x,y).
\end{equation}
Such expectations can be expressed easily by the general theory of
determinantal processes
\begin{multline*}
  \mathbb{E}\sum_{x,y}g_r(x,y)=\mathbb{E}\sum_{x\neq y}g_r(x,y)+
  \mathbb{E}\sum_{x}g_r(x,x)\\
  =\iint\limits_{(\mathbb{FP}^{d-1})^2}
  g_r(x,y)\left(K_N(x,x)K_N(y,y)-K_N(x,y)K_N(y,x)\right)\,d\sigma(x)\,d\sigma(y)
  \\+K_N(\cdot,\cdot)g_r(\cdot,\cdot).
\end{multline*}
Since all functions occurring in this last equation are zonal, we can use
\eqref{eq:zonal} to further simplify this equation
\begin{multline*}
  \mathbb{V}(\mathscr{X}_N^H,r)=\\C_{\alpha,\beta}
  \int_0^{\frac\pi2}g_r(\cos(2\theta))
  \left(K_N(1)^2-K_N(\cos(2\theta))^2\right)
  \sin(\theta)^{2\alpha+1}\cos(\theta)^{2\beta+1}\,d\theta\\+K_N(1)g_r(1);
\end{multline*}
here we have made use of the convention to write $f(\cos(2\theta(x,y)))=f(x,y)$
for zonal functions.
Using the facts that
\begin{align*}
  &C_{\alpha,\beta}\int_0^{\frac\pi2}K_N(\cos(2\theta))^2
  \sin(\theta)^{2\alpha+1}\cos(\theta)^{2\beta+1}\,d\theta=K_N(1)\\
  &C_{\alpha,\beta}\int_0^{\frac\pi2}g_r(\cos(2\theta))
  \sin(\theta)^{2\alpha+1}\cos(\theta)^{2\beta+1}\,d\theta=0
\end{align*}
we obtain
\begin{multline*}
  \mathbb{V}(\mathscr{X}_N^H,r)=\\C_{\alpha,\beta}\int_0^{\frac\pi2}
  \left(g_r(1)-g_r(\cos(2\theta))\right)K_N(\cos(2\theta))^2
  \sin(\theta)^{2\alpha+1}\cos(\theta)^{2\beta+1}\,d\theta.
\end{multline*}

We now remark that
\begin{multline*}
  g_r(1)-g_r(\cos(2\theta(x,y)))=\frac12\sigma(B(x,r)\triangle B(y,r))\\\leq
  C'\theta(x,y)\mathrm{surface}(\partial B(\cdot,r))\leq C''\sin(\theta)r^{D-1}
\end{multline*}
using \eqref{eq:symmetric}. Thus we have
\begin{equation*}
  \mathbb{V}(\mathscr{X}_N^H,r)=\mathcal{O}\left(r^{D-1}\int_0^{\frac\pi2}
  K_N(\cos(2\theta))^2
  \sin(\theta)^{2\alpha+2}\cos(\theta)^{2\beta+1}\,d\theta\right),
\end{equation*}
and it remains to estimate the integral
\begin{equation*}
  \frac{(\alpha+\beta+2)_N^2}{(\beta+1)_N^2}\int_0^{\frac\pi2}
  P_N^{(\alpha+1,\beta)}(\cos(2\theta))^2\sin(\theta)^{2\alpha+2}
  \cos(\theta)^{2\beta+1}\,d\theta;
\end{equation*}
notice that this is an integral over the square of a Jacobi polynomial against
the ``wrong'' weight function. We use the estimate \eqref{eq:jacobi-bound} for
$a=\alpha+1$ and $b=\beta$ in
the interval $[\frac cN,\frac\pi2]$ and the trivial bound
$|P_N^{(\alpha+1,\beta)}(\cos(2\theta))|\leq P_N^{(\alpha+1,\beta)}(1)$ in the
interval $[0,\frac CN]$ to obtain
\begin{multline*}
  \int_0^{\frac\pi2}
  P_N^{(\alpha+1,\beta)}(\cos(2\theta))^2\sin(\theta)^{2\alpha+2}
  \cos(\theta)^{2\beta+1}\,d\theta\\=
  \mathcal{O}\left(\int_0^{\frac CN}N^{2\alpha+2}\theta^{2\alpha+2}
    \,d\theta\right)+
\mathcal{O}\left(\int_{\frac
    CN}^{\frac\pi2}\frac1{N\sin(\theta)}\,d\theta\right)\\=
\mathcal{O}\left(\frac1N\right)+\mathcal{O}\left(\frac{\log(N)}N\right)=
\mathcal{O}\left(\frac{\log(N)}N\right).
\end{multline*}
Summing up, we have obtained
\begin{equation*}
  \mathbb{V}(\mathscr{X}_N^H,r)=
  \mathcal{O}\left(r^{D-1}\frac{(\alpha+\beta+2)_N^2}{(\beta+1)_N^2}
    \frac{\log(N)}N\right)=\mathcal{O}\left(r^{D-1}\frac{K_N(1)\log(N)}N\right).
\end{equation*}
Thus we have shown the following theorem.
\begin{theorem}
  The harmonic process $\mathscr{X}_N^H$ is hyperuniform for large and small
  balls. For balls at threshold order the weaker relation
  \begin{equation*}
    \mathbb{V}(\mathscr{X}_N^H,rK_N(1)^{-\frac1D})
    =\mathcal{O}\left(r^{D-1}\log(N)\right)
  \end{equation*}
  holds.
\end{theorem}
%%% Local Variables:
%%% mode: latex
%%% TeX-master: "Main"
%%% End: